\documentclass[a4paper,11pt]{article}
\usepackage{graphicx}
\usepackage{amssymb}
\usepackage{amsmath}
\numberwithin{equation}{section}

\begin{document}
\renewcommand{\theequation}{\thesection.\arabic{equation}}
\newcommand{\be}{\begin{equation}}
\newcommand{\ee}{\end{equation}}

\begin{center}
\Large\bf { Regularity  and geometric  character of solution of
a degenerate  parabolic equation }\\

\normalsize ~~\\
  Jiaqing Pan
\footnote {Email:  jqpan@jmu.edu.cn.   }
 \\
Institute of Mathematics, Jimei University, \\
Xiamen, 361021, P.R.China \\
 \end{center}

\vspace{0.5 cm} {\bf Abstract:}    This work studies the
regularity  and the geometric significance of solution of the
Cauchy problem for a
    degenerate  parabolic equation $u_{t}=\Delta{}u^{m}$.
 Our main objective is to improve the H$\ddot{o}$lder  estimate
 obtained by pioneers 
 and then, to show the geometric  characteristic of  free
 boundary of degenerate parabolic equation.
  To be exact, the present work will
 show that:

(1) the weak solution
$u(x,t)\in{}C^{\alpha,\frac{\alpha}{2}}(\mathbb{R}^{n}\times\mathbb{R}^{+})$,
where  $\alpha\in(0,1)$ when $m\geq2$ and $\alpha=1$ when
$m\in(1,2)$;

 (2)  the surface $\phi=(u(x,t))^{\beta}$ is a  complete   Riemannian manifold, which is tangent to $\mathbb{R}^{n}$ at
the  boundary of  the positivity set of $u(x,t)$.

 (3) the function $(u(x,t))^{\beta}$  is a classical solution
 to  another degenerate parabolic equation if
 $ \beta$ is large sufficiently;

 Moreover, some
explicit expressions about the speed of propagation of $u(x,t)$
and the continuous dependence on the nonlinearity of the equation
 are obtained.\\

Recalling the older H$\ddot{o}$lder estimate
($u(x,t)\in{}C^{\alpha,\frac{\alpha}{2}}(\mathbb{R}^{n}\times\mathbb{R}^{+})$
 with
 $0<\alpha<1$  for all $m>1$), we see our result (1) improves the
 older result and, based on this conclusion,  we can obtain (2), which shows the geometric  characteristic of  free
 boundary.
\\
~\\

  \\

{\bf AMS(2010) Subject Classifications:
}~35K15,~35K55,~35K65,~53C25

\section{Introduction}
\setcounter{equation}{0} ~~~~~~~~~

Consider   the Cauchy problem of nonlinear  parabolic equation
\begin{eqnarray}
\left\{
\begin{array}{ll}
 u_{t}=\Delta{}u^{m}~~~~~~~~~~~~~~~~~~~~~~~~~\mbox{in}~Q,\\
 u(x,0)=u_{0}(x)  ~~~~~~~~~~~~~~~~~~\mbox{on}~\mathbb{R}^{n},
 \end{array}
\right.
\end{eqnarray}
where  $Q = \mathbb{R}^{n}\times\mathbb{R}^{+},n\geq1, m>1 $ and
\begin{eqnarray} 0\leq u_{0}(x)\leq{}M,~~~~~~~~~~~~ 0 <
\int_{\mathbb{R}^{n}} u_{0}(x)dx <\infty.
 \end{eqnarray}
 The equation in (1.1) is an example of nonlinear
evolution equations and
  many interesting results, such as the existence, uniqueness,  continuous dependence on the nonlinearity of the equation
   and     large time behavior  are obtained  during the past
several decades. By a weak solution of (1.1), (1.2) in $Q$, we
mean a nonnegative  function $u(x,t)$ such that, for any given
$T>0$,
$$
\int_{0}^{T}\int_{\mathbb{R}^{n}}\left(u^{2}+|\nabla{}u^{m}|^{2}\right)dxdt<\infty
$$
and
$$
\int_{0}^{T}\int_{\mathbb{R}^{n}}\left(\nabla{}u^{m}\cdot\nabla{}f-uf_{t}\right)dxdt=
\int_{\mathbb{R}^{n}}u_{0}(x)f(x,0)dx
$$
for any continuously differentiable function $f(x,t)$ with compact
support in $ \mathbb{R}^{n}\times(0,T) $.

We know that (see  \cite{cv,arosn,
arosn2,arosn3,arosen4,pjq,jlv,jlv3}) the Cauchy problem
(1.1),(1.2) permits a unique weak solution
 $u(x,t)$ which has the following properties:
\begin{eqnarray}
&&0\leq{}u(x,t)\leq{}M,\\
&&\int_{\mathbb{R}^{n}}u(x,t)dx=\int_{\mathbb{R}^{n}}u_{0}(x)dx,\\
&&\frac{\partial{}u}{\partial{}t}\geq\frac{-u}{(m-1)t},\\
&&\Delta{}\left(\frac{m}{m-1}u^{m-1}\right)\geq\frac{-n}{n(m-1)+2}\cdot\frac{1}{t}\\
&&\|u-v\|_{L^{2}(Q_{T})}\leq{}C\max_{s\in[0,M]}\left|s^{\frac{1}{j}}-s^{\frac{m}{j}}\right|,
 \end{eqnarray}
 where $j=1,2,3...$ and $v$ is the solution to the  Cauchy problem of linear heat
 equation with the same initial value
\begin{eqnarray}
\left\{
\begin{array}{ll}
 v_{t}=\Delta{}v~~~~~~~~~~~~~~~~~~~~~~~~~~~\mbox{in}~Q,\\
 v(x,0)=u_{0}(x)
 ~~~~~~~~~~~~~~~~~~\mbox{on}~\mathbb{R}^{n}.
 \end{array}
\right.
\end{eqnarray}
Moreover,
   the solution $u(x,t)$
 can be obtained (see \cite{es,caffa}) as a limit of solutions
$u_{\eta}$($\eta\longrightarrow0^{+}$) of the Cauchy problem
\begin{eqnarray}
\left\{
\begin{array}{ll}
 u_{t}=\Delta{}u^{m}~~~~~~~~~~~~~~~~~~~~~~~~~~~~~~\mbox{in}~Q,\\
 u(x,0)=u_{0}(x)+\eta
 ~~~~~~~~~~~~~~~~~~\mbox{on}~\mathbb{R}^{n},
 \end{array}
\right.
\end{eqnarray}
and the solutions $u_{\eta}$ is taken in the classical sense. We
know that D.G.Aronson, Ph.Benilan (see
 theorem 2, p.104 in \cite{arosen4}) claimed that: if $u$
is the weak solution to the Cauchy problem (1.1) with the initial
value (1.2), then $u\in{}C(Q)$ and $u\geq0$; J. L. Vazquez  (see
Proposition 6 in Ch.2 of \cite{jlv}) proved
$u\in{}C^{\infty}(Q_{+})$, where
$$
 Q_{+}=\{(x,t)\in\mathbb{R}^{n}\times\mathbb{R}^{+}:~u(x,t)>0\}.$$
Before this,  the same conclusion was established by A. Friedman
(see theorem 11 and  corollary 2 in Ch.3 \cite{Freid} ). Moreover,
employing so called $bootstrap~argument$, D.G.Aronson, B.H.Gilding
and  L. A. Peletier (see \cite{BHG1,arosn,arosn2,arosn3}) also
claimed $u\in{}C^{\infty}(Q_{+})$ with more details. Therefore,
 we can divide the space-time
$\mathbb{R}^{n}\times\mathbb{R}^{+}$ into two parts:
$Q=Q_{+}\cup{}Q_{0}$, where
$$
Q_{0}=\{(x,t)\in\mathbb{R}^{n}\times\mathbb{R}^{+}:~u(x,t)=0\}.
$$
  Furthermore, if $
Q_{0} $ contains an   open set, say, $Q_{1}$, we can also obtain
$u(x,t)\in{}C^{\infty}(Q_{1})$ owing to $u(x,t)\equiv0$ in
$Q_{1}$. Thereby,  we may  suspect  that the solution of
degenerate parabolic equation is  actually smooth in $Q$ except a
set of measure 0. In order to improve the regularity of $u(x,t)$,
many authors have made hard effort in this direction. The earliest
contribution to the subject was made, maybe, by D.G.
  Aronson and B.H.Gilding and L.A.Peleiter(see \cite{BHG1,
arosn}). They proved
 that the
  solution to the Cauchy problem
\begin{eqnarray}
\left\{
\begin{array}{ll}
  \frac{\partial{}u}{\partial{}t}=\frac{\partial^{2}{}u^{m}}{\partial{}x^{2}}~~~~~~~~~~~~~
  ~~~~~~~~~~~~~~\mbox{in}~\mathbb{R}^{1}\times\mathbb{R}^{+},~m>1,\\
 u(x,0)=u_{0}(x)
 ~~~~~~~~~~~~~~~~~~~~~\mbox{on}~\mathbb{R}^{1}
 \end{array}
\right.
\end{eqnarray}
 is continuous in $ \mathbb{R}^{1}\times(0,+\infty)$ if the nonnegative initial value
  satisfies a good condition. Moreover,
 if the initial
value $0\leq{}u_{0}\leq{}M$ and $u_{0}^{m}$
 is Lip-continuous,
$ u(x,t)$ can   be continuous on $
\mathbb{R}^{1}\times[0,+\infty)$ (see \cite{BHG1}). As to the case
of $n\geq1$,   L.Caffarelli and A.Friedman (see \cite{ caffa}),
 proved that the solution $u(x,t)$
 to the Cauchy problem
(1.1), (1.2)
 is H$\ddot{o}$lder continuous in
$Q$ :
  \begin{eqnarray}
~~~~~~~~~~~~~~~~~~~~~~~~~~~~~~~|u(x,t)-u(x_{0},t_{0})|\leq{}C\left(|x-x_{0}|^{\alpha}+|t-t_{0}|^{\frac{\alpha}{2}}\right)~~~~~~~~~~~~~~~~~~~~~~~~~~~~(*)\nonumber
\end{eqnarray}
for some $0<\alpha<1$ and $C>0$. Moreover, for the general
equation $u_{t}=\nabla\cdot(u\nabla{}p),~p=\kappa(u)$, L
.Caffarelli and J.L.Vazques (\cite{cl2}) established the property
of finite propagation and the persistence of positivity, where
$\kappa$ may be a general operator. To study this problem more
precisely, D. G. Aronson, S.B. Angenent and J. Graveleau (see
\cite{aros,sb}) constructed a interesting radially symmetric
solution $u(r,t)$ to the focusing problem for the equation of
(1.1).  Denoting    the  porous medium pressure
$V=\frac{m}{m-1}u^{m-1}$,   they claimed   $V=Cr^{\delta}$ at the
fusing time, where $0<\delta<1,C$ is a positive constant.
Moreover,
$$
\lim_{r\downarrow0}\frac{V(r,\eta{}r^{\alpha})}{r^{2-\alpha}}=\frac{\varphi(c^{\ast}\eta)}{-\eta}.
$$

  To study  the regularity
 of the weak solution of (1.1), (1.2),
 the present work will show the following  more   precise
 conclusion: for every $h\in(m-1,m)$, there exists  a $C>0$ such that
  \begin{eqnarray}
|u(x,t)-u(x_{0},t_{0})|\leq{}C\left(|x-x_{0}|^{\frac{1}{h}}+|t-t_{0}|^{\frac{1}{2h}}\right),
\end{eqnarray}
where,
\begin{eqnarray}
~~~~~~~~~~~~h=\left\{
\begin{array}{ll}
1~~~~~~~~~~~~~~~~~~~~~~~~~~~~~~~~~~\mbox{if}~~~1<m<2,\\
 h\in(m-1,m)  ~~~~~~~~~~~~~~~~\mbox{if}~~~m\geq2.
 \end{array}
\right.\nonumber
\end{eqnarray}
We see that the range of $\frac{1}{h}$ is  $(0,1]$ not $(0,1)$,
thereby, the older  H$\ddot{o}$lder estimate $(*)$ is improved by
(1.11).

Moreover,  we will show that the functions
$\frac{\partial{u^{\beta}}}{\partial{}x_{i}} $ are continuous
 if $\beta$ is large
sufficiently because we can employ (1.11) to obtain
$$|u^{\beta}(x,t)-u^{\beta}(x_{\ast},t_{\ast})|\leq{}C \left(|x-x_{\ast}|^{\frac{\beta}{h}}+|t-t_{\ast}|^{\frac{\beta}{2h}}\right)$$
for    another positive constant $C$. By this inequality,  we want
to get the continuous partial derivatives
$\frac{\partial{}u^{\beta}}{\partial{}t}$ and
$\frac{\partial{}u^{\beta}}{\partial{}x_{i}}~i=1,2,...,n$, that is
to say,
\begin{eqnarray}
\phi(x,t)\in{}C^{1}(\mathbb{R}^{n})
\end{eqnarray}
for every given $t>0$, where $\phi(x,t)=u^{\beta}(x,t)$. In
particular, we will prove  that  the function $\phi(x,t)$
satisfies the degenerate parabolic equation
 \begin{eqnarray}
 \frac{\partial\phi_{}}{\partial{}t}=m\left[\phi^{\frac{m-1}{\beta}}\Delta{}\phi+
 \frac{m-\beta}{\beta}\phi^{\frac{m-\beta-1}{\beta}}|\nabla\phi|^{2}\right]\nonumber
\end{eqnarray}
in the classical sense.
\\

For every fixed  $t>0$, we define a n-dimensional  surface $S(t)$,
which  floats in the
 space $ \mathbb{R}^{n+1} $ with the time $t$:
\begin{eqnarray}
~~~~~~~~~~S(t):\left \{
\begin{array}{ll}
 x_{i}=x_{i},~~~~~~~~~~~~~i=1,2,3,...,n,~\\
x_{n+1}=\phi(x,t),\nonumber
 \end{array}
\right.
\end{eqnarray}
where the function $\phi(x,t)$ is mentioned above. We will discuss
the geometric  character of $S(t)$. We know that Y.Giga and
R.V.Kohn studied the  fourth-order total variation flow and the
fourth-order surface diffusion law
 (see \cite{yg}), and proved that the solution becomes identically
 zero in finite time. To be exact,   the solution  surface
 will coincide with $\mathbb{R}^{n}$ in finite time. Because this phenomenon
 will never occur for our surface $S(t)$ owing to
 (1.2) and (1.4), so we will discuss the relationship between
 $S(t)$ and $\mathbb{R}^{n}$.

Let
\begin{eqnarray}
 \left \{
\begin{array}{ll}
 g_{1}=(1,0,...,\frac{\partial{}\phi}{\partial{}x_{1}}),\\
g_{2}=(0,1,...,\frac{\partial{}\phi}{\partial{}x_{2}}),\\
 ......,\\
  g_{n}=(0,0,...1,\frac{\partial{}\phi}{\partial{}x_{n}}).
 \end{array}
\right.
\end{eqnarray}
 Define the  Riemannian metric
  on $S(t)$:
$$
(ds)^{2}=\sum\limits_{i,j=1}^{n}g_{ij}dx_{i}dx_{j} ,$$  where
$g_{ij}={g}_{i}\cdot{g}_{j}$. Clearly,
\begin{eqnarray}
 (ds)^{2}&
 =&\sum_{i=1}^{n}(1+\phi^{2}_{x_{i}})(dx_{i})^{2}+\sum_{i\neq{}j,i,j=1}^{n}\phi_{x_{i}}\phi_{x_{j}}
 dx_{i}dx_{j}\nonumber\\
 &=&\sum_{i=1}^{n}(dx_{i})^{2}+\left(\sum_{i=1}^{n}\phi_{x_{i}}dx_{i}\right)^{2}.\nonumber
\end{eqnarray}
Recalling  $
\sum\limits_{i=1}^{n}\phi_{x_{i}}dx_{i}=d\phi=dx_{n+1} $ for fixed
$t>0$, we get
\begin{eqnarray}
 (ds)^{2}
 =\sum_{i=1}^{n+1}(dx_{i})^{2}=   \sum_{i=1}^{n}(dx_{i})^{2} +(d\phi) ^{2}
 .\nonumber
\end{eqnarray}
If the derivatives $\frac{\partial{}\phi}{\partial{}x_{i}}$ are
bounded for $i=1,2,...,n$, then we can get a positive constant
$C$, such that
$\sum\limits_{i=1}^{n}\left(\phi_{x_{i}}dx_{i}\right)^{2}\leq{}C\sum\limits_{i=1}^{n}\left(dx_{i}\right)^{2}$.
Denoting  $(d\rho)^{2}=\sum\limits_{i=1}^{n}(dx_{i})^{2}$, which
is the Euclidean metric on $\mathbb{R}^{n}$,  we get
\begin{eqnarray}
  (d \rho)^{2}\leq
 (ds)^{2}
 \leq(1+C)(d\rho)^{2}.
\end{eqnarray}
 As a consequence of (1.14), we see that
the completeness of $\mathbb{R}^{n}$ yields the completeness of
$S(t)$ and therefore,
 $S(t)$  is
 a complete  Riemannian manifold. On the other hand, if we can
 obtain
\begin{eqnarray}
\nabla{}\left(\phi(x,t)\right)|_{\partial{}H_{u}(t)}=0
\end{eqnarray}
for every fixed  $t>0$, where $H_{u}(t)$ is the positivity set of
 $u(x,t)$:
$$~~~~~~~~~~~~~~~~~
H_{u}(t)=\{x\in\mathbb{R}^{n}:~u(x,t)>0\}~~~~~~~~~~~~~~~~~~~~t>0,
$$
then (1.15) encourage us to prove that:  the manifold $S(t)$ is
tangent to $\mathbb{R}^{n}$   on
 $\partial{}H_{u}(t)$.\\

It is well-known that the function
$$
v(x,t)=\left(\frac{1}{2\sqrt{\pi{}t}}\right)^{-n}\int_{\mathbb{R}^{n}}u_{0}(\xi)e^{-\frac{(x-\xi)^{2}}{4t}}d\xi
$$
 is  the solution of the Cauchy problem of the linear
heat equation (1.8) and
 $v(x,t)>0$ in $Q$  everywhere if only the initial value $u_{0}$
satisfies (1.2). This fact shows that the speed of propagation of
$v(x,t)$ is infinite, that is to say,
\begin{eqnarray}
\sup_{x\in{}H_{v}(t)}|x|=\infty.
 \end{eqnarray}
However, the degeneracy of the equation in (1.1) causes an
important phenomenon to occur, i.e. finite speed of propagation of
disturbance.  We have observed this phenomenon on the $source-
type$ solution $B(x,t;C)$ (see \cite{Ba}), where
\begin{eqnarray}
B(x,t,C)=t^{-\lambda}\left(C-\kappa\frac{|x|^{2}}{t^{2\mu}}\right)_{+}^{\frac{1}{m-1}}
 \end{eqnarray}
is the equation in (1.1) with a initial mass $M\delta(x)$, and
$$
\lambda=\frac{n}{n(m-1)+2},~~~~~~~~~~\mu=\frac{\lambda}{n},~~~~~~~~~~~~\kappa=\frac{\lambda(m-1)}{2mn}.
$$
We see that  the function $B(x,t;C)$ has compact support in space
for every fixed time. More precisely,  if $u(x,t)$ is the solution
of (1.1), (1.2), then
\begin{eqnarray}
\sup_{x\in{}H_{u}(t)}|x|={}O\left(t^{\frac{1}{n(m-1)+2}}\right)
\end{eqnarray}
when $t$ is large enough (see Proposition 17 in \cite{jlv}).
Comparing (1.16) and (1.18) and recalling the mass conservation
$\int_{\mathbb{R}^{n}}u(x,t)dx
=\int_{\mathbb{R}^{n}}u_{0}(x)dx=\int_{\mathbb{R}^{n}}v(x,t)dx$,
 we will
prove that  the solution continuously depends  on the
nonlinearity of the equation (1.1): $$\|u(\cdot,t)-v(\cdot,t)\|^{2}_{L^{2}(|x|\leq{}k)}\leq{}C\left[(m-1)+\frac{1}{k}\right].$$\\

We read the main conclusions of the present work as follows:\\

{\bf Theorem 1} {\it Assume  $u(x,t)$ be the  solution to (1.1),
(1.2).
Then $u(x,t)\in{}C(Q)$ and\\

(1) for  every given $\tau>0,K>0$ , there exists a positive $\nu$
such that
$$|u(x_{1},t_{1})-u(x_{2},t_{2})|\leq{}\nu\left(|x_{1}-x_{2}|^{\frac{1}{h}}+|t_{1}-t_{2}|^{\frac{1}{2h}}\right)$$
where $  |x_{i}|\leq{}K,t_{i}\geq\tau,~i=1,2$,
\begin{eqnarray}
~~~~~~~~~~~~h=\left\{
\begin{array}{ll}
1~~~~~~~~~~~~~~~~~~~~~~~~~~~~~~~~~~\mbox{if}~~~1<m<2,\\
 h\in(m-1,m)  ~~~~~~~~~~~~~~~~~\mbox{if}~~~m\geq2;
 \end{array}
\right.
\end{eqnarray}

 (2)  for every  $\beta>h$, the
 function
 $\phi=(u(x,t))^{\beta}\in{}C^{1}(Q)$
 and the surface $\phi=\phi(x,t)$ is a
complete Riemannian-manifold which is tangent to $\mathbb{R}^{n}$
on
$\partial{}H_{u}(t)$ for every fixed $t>0$, $h$ is defined by (1.19);\\

(3) if $\beta>2h$,  the function $\phi(x,t)$ satisfies   the
degenerate parabolic equation
  \begin{eqnarray}
 \frac{\partial\phi_{}}{\partial{}t}=m\left[\phi^{\frac{m-1}{\beta}}\Delta{}\phi+
 \frac{m-\beta}{\beta}\phi^{\frac{m-\beta-1}{\beta}}|\nabla\phi|^{2}\right]\nonumber
\end{eqnarray}
in the classical sense in $Q$.
}\\

{\bf Theorem 2} {\it Assume  $u(x,t)$ be the  solution to (1.1),
(1.2), $B_{\delta}=\{x\in\mathbb{R}^{n}:~|x|<\delta\}$ for some
$\delta>0$.
 If supp $u_{0}\subset{}B_{\delta}$, then for every given
$t>0$,
\begin{eqnarray}
\sup_{x\in{}H_{u}(t)}|x|\geq\chi(t),
\end{eqnarray}
where,
 \begin{eqnarray}
\chi(t)=\left[
  (m-1)\pi ^{\frac{(1-m)n}{2}}
\Gamma(1+\frac{n}{2})^{m-1}\cdot\left(\int_{\mathbb{R}^{n}}u_{0}dx\right)^{m-1}t\right]^{\frac{1}{2+(m-1)n}}.\nonumber
\end{eqnarray}
Moreover, for every given $T>0$, there is a positive
$C_{\ast}=C_{\ast}(T)$ such that
\begin{eqnarray}
 \int_{|x|\leq{}k} \left[ v(x,t)-u(x,t)
\right]^{2} dx
 \leq{}C_{\ast}\left[(m-1)+\frac{1}{k}\right]
\end{eqnarray}
with respect to $t\in(0,T)$ uniformly, where  $v(x,t)$ is the
solution of (1.8).
 }\\

Let us also note that the manifolds $S(t)$ and $\mathbb{R}^{n}$
are two surfaces in $\mathbb{R}^{n+1}$ and the Cauchy problem
(1.1), (1.2) can be regarded as a mapping
$\Phi(t):~\mathbb{R}^{n}\longrightarrow~S(t)$. Thus, besides the
theorems mentioned  above, we will give an example to show the
intrinsic
properties about the manifold $S(t)$.\\

\section{The proof of Theorem 1}
\setcounter{equation}{0} ~~~~~~~~~

{\bf Lemma 1} {\it If  $u(x,t)$ is the   weak solution to (1.1),
(1.2) in $Q$. Then for every  $ h\in(m-1,m)$, there is a
$C_{1}=C_{1}(h,m,M)$ such that
\begin{eqnarray}
~~~~~~~~~~~~~~~~~~~~\left|\nabla{}u^{h}\right|^{2}\leq{}\frac{1}{{C_{1}t}}~~~~~~~~~~~~~~~~~~\mbox{in}~Q
\end{eqnarray}
in the sense of distributions in $Q$.}
 \\

{\bf Proof:} We first prove (2.1) for the classical solutions
$u_{\eta}(x,t)$. Set
$$~~~~~~~~~~~~~~~~~~~~~~~~~u_{\eta}^{m}=V^{q}~~~~~~~~~~~~~\mbox{for}~ q\in\left(1,\frac{m}{m-1}\right).$$
Then
$$
V_{t}=mV^{q-\frac{q}{m}}\Delta{}V+m(q-1)V^{q-1-
\frac{q}{m}}|\nabla{}V|^{2}.
$$
 Differentiating  this equation with respect to $x_{j}$ and
 multiplying
though by $\frac{\partial{}V}{\partial{}x_{j}}$, letting
$h_{j}=\frac{\partial{}V}{\partial{}x_{j}}$, we get
\begin{eqnarray}
\frac{1}{2}(h_{j}^{2})_{t}&=&mV^{q-\frac{q}{m}}h_{j}\Delta{}h_{j}+
m(q-\frac{q}{m})V^{q-1-\frac{q}{m}}h_{j}^{2}\Delta{}V\nonumber\\
&&+m(q-1)(q-1-\frac{q}{m})V^{q-2-\frac{q}{m}}h_{j}^{2}|\nabla{}V|^{2}+2m(q-1)V^{q-1-\frac{q}{m}}h_{j}\nabla{}V\cdot\nabla{}h_{j}\nonumber
\end{eqnarray}
for $ j=1,2,...,n  $. Setting $$H^{2}=\sum_{j=1}^{n}h_{j}^{2}~,$$
we obtain
\begin{eqnarray}
&&H^{2}_{t}-mV^{q-\frac{q}{m}}\Delta{}H^{2}\nonumber\\
&=& 2m(q-\frac{q}{m})V^{q-1-\frac{q}{m}}H^{2}\Delta{}V
+2m(q-1)(q-1-\frac{q}{m})V^{q-2-\frac{q}{m}}\sum_{j=1}^{n}h_{j}^{2}|\nabla{}V|^{2}\nonumber\\
&&+
2m(q-1)V^{q-1-\frac{q}{m}}\nabla{}V\cdot\sum_{j=1}^{n}\nabla{}h^{2}_{j}-2mV^{q-\frac{q}{m}}\sum_{j=1}^{n}|\nabla{}h_{j}|^{2}\nonumber\\
&\leq& 2
  m(q-\frac{q}{m})V^{q-1-\frac{q}{m}}H^{2}\Delta{}V
  +2m(q-1)(q-1-\frac{q}{m})V^{q-2-\frac{q}{m}}H^{4}\nonumber\\
  && +2
   m(q-1)V^{q-1-\frac{q}{m}}\nabla{}V\cdot\nabla{}H^{2}.\nonumber
\end{eqnarray}
It follows from $ q\in\left(1,\frac{m}{m-1}\right)$ that
$(q-1)(q-1-\frac{q}{m})<0$ and ${q-2-\frac{q}{m}}<0$. Hence,
\begin{eqnarray}
2m(q-1)(q-1-\frac{q}{m})V^{q-2-\frac{q}{m}}&\leq&{}-\widetilde{C}_{1},\nonumber
\end{eqnarray}
where $
\widetilde{C}_{1}=2m\left|(q-1)(q-1-\frac{q}{m})\right|(M+\eta)^{m-\frac{2m}{q}-1}.
$
Setting
\begin{eqnarray} L({H}^{2})&=&
mV^{q-\frac{q}{m}}\Delta{}{H}^{2} -\widetilde{C}_{1}{H}^{4}+
2m(q-\frac{q}{m})V^{q-1-\frac{q}{m}}{H}^{2}\Delta{}V\nonumber\\
&&+ 2m(q-1)V^{q-1-\frac{q}{m}}\nabla{}V\cdot\nabla{}{H}^{2},
\nonumber
\end{eqnarray}
  we get
$$
H^{2}_{t}\leq{}L(H^{2}).
$$
Taking
$h^{\ast}_{1}=(\frac{1}{\widetilde{C}_{1}t})^{\frac{1}{2}},h^{\ast}_{i}=0$
for $i=2,3,...,n$, and setting $
Q_{\ast}^{2}=\sum\limits_{i=1}^{n}h_{i}^{\ast}=
\frac{1}{\widetilde{C}_{1}t}$, we see that
 the function $ Q_{\ast}^{2}$ is a solution to the
equation
\begin{eqnarray}
\frac{\partial}{\partial{}t}Q_{\ast}^{2}={}L(Q_{\ast}^{2})\nonumber
\end{eqnarray}
with the initial condition $Q^{2}_{\ast}(0)=+\infty$.

On the other hand, let $u_{\eta}^{m}=K$,
 then $K$ is the solution to the  Cauchy problem of the linear
 parabolic equation
\begin{eqnarray}
~~~~~~~~~~~~~~~~~~~~~~~~~~~~~~~~~~~~~~K_{t}=g(x,t)\Delta{}K~~~~~~~~~~~~~~~~~~~~~~~~~\mbox{in}~Q\nonumber
\end{eqnarray}
with the inial value $(u_{0}+\eta)^{m}$, where,
$g(x,t)=m(u_{\eta}(x,t))^{m-1}$ and $u_{\eta}(x,t)$ is a known
function, which is the classical
 solution to (1.9) and $\eta\leq{}u_{\eta}\leq{}M+\eta$.
By Theorem 5.1 in \cite{OA}, we get
$K\in{}H^{2+\alpha,1+\frac{\alpha}{2}}(\mathbb{R}^{n}\times(0,T))$
for any $T>0$ (Even if $u_{0}$ does not have the required
smoothness we may approximate it ( by mollification ) with smooth
functions $u_{0\eta}$). Therefore, $\Delta{}K,
\frac{\partial{}K}{\partial{}t}$ and
$\frac{\partial{}K}{\partial{}x_{i}} (i=1,2,...,n)$ are bounded.
To be exact, there is a positive $\mu_{0}$, which may depend on
$\eta$, such that
$$
|\Delta{}K|+|\frac{\partial{}K}{\partial{}t}|+|\nabla{}K|\leq{}\mu_{0}.
$$
Therefore,
\begin{eqnarray}
|\nabla{}V|&=&|\nabla{}u_{\eta}^{\frac{m}{q}}|\nonumber\\
&=&|\frac{1}{q}u_{\eta}^{\frac{m}{q}-m}\nabla{}K|\nonumber\\
&\leq&{}\frac{1}{q}\eta^{(\frac{1}{q}-1)m}\mu_{0}.\nonumber\end{eqnarray}
Letting
 $C'=\frac{1}{q}\eta^{(\frac{1}{q}-1)m}\mu_{0}$, we obtain $|\nabla{}V|\leq{}C'$.
Similarly, we can get positive $C''$, which may depend on $\eta$
also, such that
 $|\Delta{}V|\leq{}C''$.

Now we can employing
  the
comparison theorem and get
\begin{eqnarray}
H^{2}\leq{}\frac{1}{\widetilde{C}_{1}t} .
\end{eqnarray}
Letting $\eta\longrightarrow0$ in (2.2) gives
$|\nabla{}u^{\frac{m}{q}}|^{2}\leq\frac{1}{{C}_{1}t}$ for $
q\in\left(1,\frac{m}{m-1}\right)$ with
$$C_{1}=2m\left|(q-1)(q-1-\frac{q}{m})\right|M^{m-\frac{2m}{q}-1}.$$
Setting $h=\frac{m}{q}$ in (2.2) yields $h\in(m-1,m)$, and (2.1)
follows. ~~~~
~~~~~~~~~~~~~~~~~~~~~~~~~~~~~~~~~~~~$\square$\\

To prove Theorem 1, we need to show an   ordinary inequality
firstly:
 \begin{eqnarray}
~~~~~~~~~~~~~~~~~|a-b|^{\beta}\leq|a^{\beta}-b^{\beta}|~~~~~~~\mbox{for}~a,b\geq0,~\beta{}>1.
\end{eqnarray}
In fact, (2.3) is right for $a=b$.  If $a>b$, we can easily get
the following inequalities:
$$
\left(1-\frac{b}{a}\right)^{\beta}\leq1-\frac{b}{a}~~~~~~~\mbox{and}~~~~~~~1-\left(\frac{b}{a}\right)^{\beta}\geq1-\frac{b}{a}
$$
thanks to $0\leq\frac{b}{a}<1$. Thereby, $
\left(1-\frac{b}{a}\right)^{\beta}\leq1-\left(\frac{b}{a}\right)^{\beta}.\nonumber
$ This inequality gives
\begin{eqnarray}
|a-b|^{\beta}=a^{\beta}\left|1-\frac{b}{a}\right|^{\beta}
\leq{}a^{\beta}-b^{\beta}.\nonumber
\end{eqnarray}
So (2.3) holds for $a>b\geq0$. Certainly, (2.3) is also right when
$0\leq{}a<b$.
\\

We are now in a position to establish our Theorem 1.\\

{\bf To prove  (1) of Theorem 1~} It follows from (2.1) that
\begin{eqnarray}
\left|u^{h}(x_{1},t)-u^{h}(x_{2},t)\right|\leq{}
(C_{1}t)^{-\frac{1}{2}}|x_{1}-x_{2}|
\end{eqnarray}
for every $(x_{1},t),(x_{2},t)\in{}Q,h\in(m-1,m)$. If $1<m<2$, we
take $h=1$, and therefore, (2.4) yields
\begin{eqnarray}
\left|u(x_{1},t)-u(x_{2},t)\right|\leq{}
(C_{1}t)^{-\frac{1}{2}}|x_{1}-x_{2}|.\nonumber
\end{eqnarray}
If $m\geq2$, we take $h>1$ for every $h\in(m-1,m)$. In this case,
we use (2.3) in (2.4) and obtain
\begin{eqnarray}
\left|u(x_{1},t)-u(x_{2},t)\right|\leq{}(C_{1}t)^{-\frac{1}{2h}}|x_{1}-x_{2}|^{\frac{1}{h}}.\nonumber
\end{eqnarray}
 Therefore, for every given  $m>1$, there always exists a suitable
 positive number
$h\geq1$ such that
\begin{eqnarray}
\left|u(x_{1},t)-u(x_{2},t)\right|\leq (C_{1}\tau)^{-\frac{1}{2h}}
|x_{1}-x_{2}|^{\frac{1}{h}}
\end{eqnarray}
for every
$(x_{1},t),(x_{2},t)\in\mathbb{R}^{n}\times[\tau,\infty)$ with
 given $\tau>0$,
 where
 $
{h}=1 ~\mbox{for}~1<m<2; ~h\in(m-1,m)
 $ for $m\geq2$.
 Employing  the
well-known theorem on the H$\ddot{o}$lder continuity with respect
to the time variable (see  \cite{BHG}), we   obtain
\begin{eqnarray}
\left|u(x,t_{1})-u(x,t_{2})\right|\leq{}\mu{}|t_{1}-t_{2}|^{\frac{1}{2h}}
\end{eqnarray}
for $|x|\leq{}K,~t_{1},t_{2}>\tau$ and $|t_{1}-t_{2}|$ small
sufficiently, where $\mu$ depends on $(C_{1}\tau)^{-\frac{1}{2h}},
K$ is any fixed positive constant. Combining (2.5) and (2.6) we
get another positive constant $\nu$ such that
\begin{eqnarray}
\left|u(x_{1},t_{1})-u(x_{2},t_{2})\right|\leq{}
\nu\left[|x_{1}-x_{2}|^{\frac{1}{h}}+|t_{1}-t_{2}|^{\frac{1}{2h}}\right]
\end{eqnarray}
for all $|x_{i}|\leq{}K,t_{i}\geq\tau,~i=1,2$, where
\begin{eqnarray}
~~~~~~~~~~~~h=\left\{
\begin{array}{ll}
1~~~~~~~~~~~~~~~~~~~~~~~~~~~~~~~~~~\mbox{if}~~~1<m<2,\\
 h\in(m-1,m)  ~~~~~~~~~~~~~~~~\mbox{if}~~~m\geq2.
 \end{array}\nonumber
\right.
\end{eqnarray}
 Certainly, this gives $u\in{}C(Q).$
~~~~~~~~~~~~~~~~~~~~~~~~~~~~~~~~~~~~~~~~~~~ ~~~$\square$\\

 {\bf To prove
(2) of Theorem 1}  Denote $\beta=h+\varepsilon$ and set
\begin{eqnarray}\phi(x,t)=u^{\beta}\nonumber
\end{eqnarray}
for every $\varepsilon>0$. It follows from (2.1) that
\begin{eqnarray}
~~~~~~~~~~~~~~~~~~~~~|\nabla{}\phi(x,t)|\leq{}C_{2}u^{\varepsilon}(x,t)t^{-\frac{1}{2}}~~~~~
~~~~~~~~~~~~~~~\mbox{in}~Q
\end{eqnarray}
for some $C_{2}>0$. To prove the function
$\phi(x,t)\in{}C^{1}(Q)$, we first see $u\in{}C^{\infty}(Q_{+})$;
second, (2.8) implies $\nabla{}\phi(x_{\ast},t_{\ast})=0$ for
every $(x_{\ast},t_{\ast})\in{}Q_{0}$. Combining (2.7) and (2.8)
gives
\begin{eqnarray}\left|\nabla{}\phi(x,t)-\nabla{}\phi(x_{\ast},t_{\ast})\right|&=&\left|\nabla{}\phi(x,t)\right|\nonumber\\
&\leq&C_{2}{}t^{-\frac{1}{2}}\nu^{\varepsilon}
\left[|x-x_{\ast}|^{\frac{1}{h}}+|t-t_{\ast}|^{\frac{1}{2h}}\right]^{\varepsilon}.
\end{eqnarray}
 This inequality tells us $\nabla\phi(x,t)\in{}C(Q)$ and therefore,
$\phi(x,t)\in{}C^{1}(Q)$. Furthermore, for every fixed $t>0$, the
continuity of $u(x,t)$ implies $H_{u}(t)$ is a open set. Thus
 $\phi(x_{\ast},t)=\nabla\phi(x_{\ast},t)=0$ for $(x_{\ast},t)\in{\partial{}H_{u}(t)}$.  This fact claims
  that the  surface $\phi=\phi(x,t)$ touches $\mathbb{R}^{n}$ at
$\partial{}H_{u}(t)$.
 In other words, $\mathbb{R}^{n}$ is just the
tangent plane of $\phi(x,t)$ at $\partial{}H_{u}(t)$.

Now we can define a n-dimensional surface $S(t)$ for every fixed
$t>0$ as before:
\begin{eqnarray}
~~~~~~~~~~S(t):\left \{
\begin{array}{ll}
 x_{i}=x_{i},~~~~~~~~~~~~~i=1,2,3,...,n,~\\
x_{n+1}=\phi(x,t).\nonumber
 \end{array}
\right.
\end{eqnarray}
 Defining the
Riemannian metric (1.13) on $S(t)$,  we get (1.14), more
precisely,
\begin{eqnarray}
  (d \rho)^{2}\leq
 (ds)^{2}
 &\leq&\left(1+\max_{i=1,2,...,n}\left|\frac{\partial\phi}{\partial{}x_{i}}\right|^{2}\right)(d\rho)^{2}.\nonumber
\end{eqnarray}
It follows from (2.8) that $
\max\limits_{i=1,2,...,n}\left|\frac{\partial\phi}{\partial{}x_{i}}\right|^{2}\leq{}\left(C_{2}
u^{\varepsilon}t^{-\frac{1}{2}}\right)^{2}$. Moreover, as  a
result of theorem 9 in Chap III of \cite{jlv}, we get
$u(x,t)\leq{}C_{3}t^{-\frac{n}{n(m-1)+2}}$ for some $C_{3}>0$.
Hence,
\begin{eqnarray}
\max\limits_{i=1,2,...,n}\left|\frac{\partial\phi}{\partial{}x_{i}}\right|^{2}\leq{}C_{4}
t^{\frac{-2n\varepsilon}{n(m-1)+2}-1}
 \end{eqnarray}
for another  positive $C_{4}$. Thus,
\begin{eqnarray}
  (d \rho)^{2}\leq
 (ds)^{2}
 &\leq&\left(1+C_{4}
t^{\frac{-2n\varepsilon}{n(m-1)+2}-1}
 \right)(d\rho)^{2}.~
\end{eqnarray}
 As a consequence of (2.11),
the completeness of $\mathbb{R}^{n}$ yields the completeness of
$S(t)$ and therefore,
 $S(t)$  is
 a complete  Riemannian manifold.
~~~~~~~~~~~~~~~~~~~~~~~~~~~~~~~~~~~~~~~~~~~~~~~~~~~~~~~~~~~~~~~~~~~$\square$
\\

{\bf To  prove (3) of Theorem 1~}  Recalling
$u(x,t)=\lim\limits_{\eta\longrightarrow0^{+}}u_{\eta}$ and
$u_{\eta}\geq\eta$, $u_{\eta}$ are the classical solutions to the
Cauchy problem (1.9), we can make
$$
\phi_{\eta}=u_{\eta}^{\beta}
$$ for $\beta=h+\varepsilon>2h$, and
then $\phi_{\eta}$ satisfies the degenerate parabolic  equation
\begin{eqnarray}
 ~~~~~~~~~~~~~~~~~~\frac{\partial\phi_{\eta}}{\partial{}t}=m\left[\phi_{\eta}^{\frac{m-1}{\beta}}\Delta{}\phi_{\eta}+
 \frac{m-\beta}{\beta}\phi_{\eta}^{\frac{m-\beta-1}{\beta}}|\nabla\phi_{\eta}|^{2}\right]~~~~~~~~~~~\mbox{in}~Q.\nonumber
\end{eqnarray}
Recalling $\lim\limits_{\eta\longrightarrow0^{+}}\phi_{\eta}=\phi$
 and
  $u\in{}C^{\infty}(Q_{+})$, we see that
  $\phi(x,t)\in{}C^{\infty}(Q_{+})$ and $\phi$ satisfies
the equation
\begin{eqnarray}
 \frac{\partial\phi_{}}{\partial{}t}=m\left[\phi^{\frac{m-1}{\beta}}\Delta{}\phi+
 \frac{m-\beta}{\beta}\phi^{\frac{m-\beta-1}{\beta}}|\nabla\phi|^{2}\right]~
\end{eqnarray}
in $Q_{+}$ . We next  prove that (2.12) is also right in $Q$. To
do this, we first see that
  (2.7) yields
\begin{eqnarray}
 \left|u(x,t)-u(x,t_{\ast})\right|=u(x,t)\leq{}\nu
 \left|t-t_{\ast}\right|^{\frac{1}{2h}}~~~~~~~~~~~~~~~~\mbox{for}~(x,t_{\ast})\in{}Q_{0}.\nonumber
\end{eqnarray}
Thus,
$|\phi(x,t)-\phi(x,t_{\ast})|=\phi(x,t)\leq{}\nu^{\beta}|t-t_{\ast}|^{\frac{\beta}{2h}}$.
Thus,
\begin{eqnarray}
 \left|\frac{\phi(x,t)-\phi(x,t_{\ast})}{t-t_{\ast}}\right|\leq\nu^{\beta}
 \left|t-t_{\ast}\right|^{\frac{\beta{}-2h}{2h}}.
\end{eqnarray}
(2.13) shows the   continuity   of the function
$\frac{\partial\phi}{\partial{}t}$, specially,
$\frac{\partial{\phi} }{\partial{}t}|_{(x,t_{\ast})}=0$ owing to
$\beta>2h$. On the other hand, similar to (2.5) and (2.8), we can
get
$$|u(x_{1},...,x_{j},...,x_{n},t)-u(x_{1},...,x_{j\ast},...,x_{n},t)|
\leq(C_{1}t)^{-\frac{1}{2h}}|x_{j}-x_{j\ast}|^{\frac{1}{h}}$$ for
every
$(x_{1},...,x_{j},...,x_{n},t),(x_{1},...,x_{j\ast},...,x_{n},t)\in{}Q$
and  $
\left|\frac{\partial{}\phi}{\partial{}x_{i}}\right|\leq{}C_{2}u^{\varepsilon}t^{-\frac{1}{2}}~\mbox{for}~i=1,2,...,n.
$ Thereby, if $ (x_{1},...,x_{j\ast},...,x_{n},t) \in{}Q_{0}$, we
have
\begin{eqnarray}
\left|\frac{\partial{}\phi(x_{1},...,x_{j},...,x_{n},t)}{\partial{}x_{i}}-\frac{\partial{}\phi(x_{1},...,x_{j_{\ast}},...,x_{n},t)}{\partial{}x_{i}}
 \right|&=&\left|\frac{\partial{}\phi(x_{1},...,x_{j},...,x_{n}t)}{\partial{}x_{i}}
 \right|\nonumber\\
 &\leq&C_{2}t^{-\frac{1}{2}}(C_{1}t)^{-\frac{\varepsilon}{2h}}|x_{j}-x_{j\ast}|^{\frac{\varepsilon}{h}}.\nonumber
\end{eqnarray}
This yields
\begin{eqnarray}
\frac{1}{|x_{j}-x_{j\ast}|}\left|\frac{\partial{}\phi(x_{1},...,x_{j},...,x_{n},t)}{\partial{}x_{i}}-\frac{\partial{}\phi(x_{1},...,x_{j_{\ast}},...,x_{n},t)}{\partial{}x_{i}}
 \right|
 \leq{}C_{2}(C_{1}t)^{-\frac{\varepsilon}{2h}}|x_{j}-x_{j\ast}|^{\frac{\varepsilon}{h}-1}\nonumber
\end{eqnarray}
for $i,j=1,2,...,n$. This gives the continuity of the function
$\frac{\partial^{2}{\phi}}{\partial{}x_{i}x_{j}}$ thanks to
$\varepsilon>h$.  In particular,
$$~~~~~~~~~~~~~~~~~~~~~~~\frac{\partial^{2}{\phi}}{\partial{}x_{i}x_{j}}=0 ~~~~~~~~
\mbox{on}~Q_{0}
$$   for $  i,j=1,2,...,n$.  It follows from $m>1$ that
$\phi^{\frac{m-1}{\beta}}\Delta{}\phi_{}=0 $ on $Q_{0}$.
Similarly,  $\phi^{\frac{m-\beta-1}{\beta}}|\nabla\phi|^{2}=0$ on
$ Q_{0}$.

Finally,
 as previously mentioned above,
  we  deduce  that the function $\phi(x,t)$ satisfies (2.12) on
  $Q_{0}$.~~
  ~~~~~~~~~~~~~~~~~~~~~~~~~~~~~~~~~~~~~~~~~~~~~~~~~~~~~~~~~~~~~~~~~~~~~~~~~~~~~~~~~~~~~~~~~~~~~~~~~~~~~~~~~~~~~~~$\square$\\

As  an applications of our Theorem 1, here we give an example to
show the large time behavior on the intrinsic properties of the
manifold $S(t)$.\\

{\bf Example~ ( the first fundamental form on $S(t)$)~ }    In
fact, the
 first fundamental form on the manifold $S(t)$ is $
(ds)^{2}=\sum\limits_{i,j=1}^{n}g_{ij}dx_{i}dx_{j} $. By (2.10),
 \begin{eqnarray}\left|(ds)^{2}-(d\rho)^{2}\right|=(d\phi)^{2}&=&\left(\sum_{i=1}^{n}{\phi_{x_{i}}}dx_{i}\right)^{2}\nonumber\\
 &\leq&C_{4}t^{-\frac{2n\varepsilon}{n(m-1)+2}-1}(d\rho)^{2}.\nonumber\end{eqnarray}
where  $ (d\rho)^{2}=\sum\limits_{i,j=1}^{n}dx_{i}^{2} $ is just
the  first fundamental form on  $\mathbb{R}^{n}$. Thus,
 \begin{eqnarray}
\frac{(ds)^{2}}{(d\rho)^{2}}=1+O\left(t^{t^{-\frac{2n\varepsilon}{n(m-1)+2}-1}
}\right)
\end{eqnarray}
when $t$ is large sufficiently.
\\

\section{The proof of Theorem 2}
\setcounter{equation}{0} ~~~~~~~~~

To prove  Theorem 2, we need to establish  a more precise
Poincar$\acute{e}$ inequality. It is well-known that if
$\lambda_{1}$ is the minimum positive eigenvalue and $\psi_{1}$ is
the corresponding eigenfunction of the Dirichlet problem
\begin{eqnarray}
\left\{
\begin{array}{ll}
 \Delta{}u=-\lambda{}u~~~~~~~~~~~~~~~~~~~ ~\mbox{in}~\Omega,\\
 u=0  ~~~~~~~~~~~~~~~~~~~~~~~~~ ~\mbox{on}~ \partial\Omega,
 \end{array}
\right.\nonumber
\end{eqnarray}
then $
{\lambda_{1}}\|\psi_{1}\|^{2}_{L^{2}(\Omega)}=\|\nabla{}\psi_{1}\|^{2}_{L^{2}(\Omega)}$,
where $\Omega$ is a boundary domain in $\mathbb{R}^{n}$.
 Moreover, if $\psi\in{}H_{0}^{1}(\Omega)$, then
 Poincar$\acute{e}$ inequality
claims that  there exists a positive constant $k$ such that
$k\|\psi\|^{2}_{L^{2}(\Omega)}\leq{}\|\nabla{}\psi\|^{2}_{L^{2}(\Omega)}$.
 According to Qiu (see p.98 in \cite{qct}), $ \lambda_{1}$ is the
maximum of all such $k$. We know that there are many kinds of
choice for such $k$. For example, Wu  (see p.13 in \cite{wzq})
proved that
\begin{eqnarray}
k\leq{}\rho^{-2}\nonumber
\end{eqnarray}
if $\Omega=
\{(x_{1},x_{2},...,x_{N})\in\mathbb{R}^{N}:~a_{i}<x_{i}<a_{i}+\rho~\}$.
In order to prove Theorem 2 we need to show that such choice is
also right if $\Omega$
is a sphere of $\mathbb{R}^{n}$.\\

Denote
$$
B=\{x\in\mathbb{R}^{n}:~|x-x_{0}|<\rho\}
$$
for  $x_{0}\in\mathbb{R}^{n}$ and  $\rho>0$. We have the following result.\\

{\bf Lemma 2 } {\it If  $\Omega\subset{}B$,  $u\in{}H^{1}(
\Omega)$ and $u(x)=0$ for $x\in\partial{\Omega}$. Then
\begin{eqnarray}~~~~~~~~~
\|u\|_{L^{2}(\Omega)}\leq{}\rho\|\nabla{}u\|_{L^{2}(\Omega)}~.
\end{eqnarray}}\\

 {\bf Proof  ~}  We first suppose $u\in{}C_{0}^{\infty}(\Omega)$.
 For every $x\in{}\Omega$, there is a
 $x_{_{\ast}}\in{}\partial{}\Omega$,
such that the three points $x_{0},~x$ and $x_{\ast}$ lie on a
radius $\overline{x_{0}~x_{\ast}}$.
 Denote the vector from  $x_{\ast}$ to $x$
by $r$. We have
\begin{eqnarray}
 u(x)&=&u(x)-u(x_{\ast})\nonumber\\
&=&\int_{x_{\star}}^{x}\frac{\partial{}u(x)}{\partial{}r}dr.~~
 \nonumber
\end{eqnarray}
Using the H$\ddot{o}$lder inequality gets
\begin{eqnarray}
 |u(x)|^{2}
&\leq&\rho\int_{x_{\ast}}^{x}\left|\frac{\partial{}u(x)}{\partial{}r}\right|^{2}dr
.
 \nonumber
\end{eqnarray}
Thus,
\begin{eqnarray}
 \int_{\Omega}|u(x)|^{2}dx
&\leq&\rho
\int_{\Omega}\int_{x_{\ast}}^{x}\left|\frac{\partial{}u(x)}{\partial{}r}\right|^{2}dr{}dx\nonumber\\
 &\leq&\rho^{2}
 \int_{\Omega}|\nabla{}u(x)|^{2}dx.
 \nonumber~
\end{eqnarray}
The general case is done by approximation.
~~~~~~~~~~~~~~~~~~~~~~~~~~~~~~~~~~~~~~~~~~~~~~~~~~~~~~~~~~~~~$
\square$
\\

{\bf To prove (1.20)~} Assume $u(x,t)$ be the solution of (1.1),
(1.2). Integrating (1.5) from $t_{1}$ to $t_{2}$ yields
$\ln{}u(x,t_{2})-\ln{}u(x,t_{1})\geq-
\frac{1}{m-1}(\ln{}t_{2}-\ln{}t_{1})$ for~$t_{2}>t_{1}$. This
means $
u(x,t_{2})\cdot{}t_{2}^{\frac{1}{m-1}}\geq{}{u(x,t_{1})}\cdot{}t_{1}^{\frac{1}{m-1}}$
 for $~t_{2}>t_{1}\geq0$.  Therefore,
\begin{eqnarray}
~\left\{
\begin{array}{ll}
 if~ u(x_{0},t_{0})=0, ~then~
u(x_{0},t)=0~~ for ~every ~~0\leq{}t<t_{0}; \\
if~ u(x_{0},t_{0})>0, ~then~ u(x_{0},t)>0 ~~for ~every ~~t>t_{0}.
 \end{array}
\right.
\end{eqnarray}
By (3.2), we see that
$$~~~~~~~~~~~H_{u}(t)\supset{}B_{\delta}~~~~~~~~~~~~~~t>0.$$
Absolutely, the proof is finished if
$\sup\limits_{x\in{}H_{u}(t)}|x|=\infty$, otherwise, we
 set
$$
K'= \gamma+\sup\limits_{x\in{}H_{u}(t)}|x|
$$
for $\gamma>0$. Thus,
\begin{eqnarray}
u(x,t)=0~~~~~~~~~~~~~~~\mbox{for}~x\in\mathbb{R}^{n}-B_{K'}.
 \end{eqnarray}
 It follows from (1.5) that
  \begin{eqnarray} \int_{ B_{K'} } u^{m}\Delta{}u^{m}dx
 \geq{}-\frac{1}{(m-1)t}
\int_{  B_{K'}}u^{1+m} dx,\nonumber
\end{eqnarray}
so that
\begin{eqnarray} \int_{ B_{K'}} |\nabla{}u^{m}|^{2}dx
 \leq{}\frac{1}{(m-1)t}
\int_{  B_{K'}}u^{1+m} dx.\nonumber
\end{eqnarray}
Using (3.1) in this inequality, we obtain
\begin{eqnarray} \int_{ B_{K'}} u^{2m}dx
 \leq{}\frac{K'^{2}}{(m-1)t}
\int_{  B_{K'}}u^{1+m} dx.\nonumber
\end{eqnarray}
Employing    the  H$\ddot{o}$lder  inverse inequality, we have
\begin{eqnarray}
\int_{B_{K'}}u^{2m}dx \geq \left(\int_{B_{K'}}u^{1+m}dx\right)^{
\frac{2m}{1+m}}|B_{K'}|^{\frac{1-m}{1+m}},\nonumber
\end{eqnarray}
where $|B_{K'}|$ is the  volume   of  $B_{K'}$ and $|B_{K'}| = \pi
^{\frac{n}{2}} \Gamma(1+\frac{n}{2})^{-1}K'^{n}$.  So we get
\begin{eqnarray}
\left(\int_{B_{K'}}u^{1+m}dx\right)^{
\frac{m-1}{1+m}}|B_{K'}|^{\frac{1-m}{1+m}}\leq\frac{K'^{2}}{(m-1)t}.
\end{eqnarray}
Using the  H$\ddot{o}$lder  inverse inequality again, we have
\begin{eqnarray}
\int_{B_{K'} }u^{1+m}dx&\geq&\left( \int_{
B_{K'}}udx\right)^{1+m}|B_{K'}|^{-m}.
\end{eqnarray}
It follows from
 (1.4) and  (3.3) that
$ \int_{\mathbb{R}^{n}}u(x,t)dx=\int_{ B_{K'}}u(x,t)dx
=\int_{\mathbb{R}^{n}}u_{0}(x)dx $. Combining (3.4) and (3.5) we
get
\begin{eqnarray}
\left(
\int_{\mathbb{R}^{n}}u_{0}dx\right)^{m-1}|B_{K'}|^{1-m}\leq\frac{K'^{2}}{(m-1)t}.
\end{eqnarray}
Now we get
\begin{eqnarray}
  (m-1)\left(\int_{\mathbb{R}^{n}}u_{0}dx\right)^{m-1}t\leq
 K'^{2+(m-1)n} \cdot\pi ^{\frac{(m-1)n}{2}}\cdot
\left(\Gamma(1+\frac{n}{2})\right)^{1-m}.\nonumber
\end{eqnarray}
Letting
 $ \gamma\longrightarrow0$ gives
\begin{eqnarray}
 ~~~~~~~~~~~~~~~~~~~~~~~~~~~~~~~~~~~~~~~~~~~~ \sup_{x\in{}H_{u}(t)} |x|\geq{} \chi(t)
  ~~~~~~~~~~~~~~t>0.~~~~~~~~~~~~~~~~~~~~~~~~~~~~~~~
 \square\nonumber
\end{eqnarray}
\\

 {\bf To prove (1.21)~} Assume  $u(x,t)$ and $v(x,t)$ be
 the
 solutions to (1.1) and (1.8) respectively.
 Employing the well-known result (see \cite{jlv}), we have
  \begin{eqnarray} \frac{1}{1+m}\int_{\mathbb{R}^{n}}u^{1+m}(x,T)dx+
  \int_{Q_{T}}|\nabla{}u^{m}|^{2}
 dxdt\leq\frac{1}{1+m}\int_{\mathbb{R}^{n}}u_{0}^{1+m}dx \end{eqnarray}
and
 \begin{eqnarray} \frac{1}{2}\int_{\mathbb{R}^{n}}v^{2}(x,T)dx+
  \int_{Q_{T}}|\nabla{}v|^{2}
 dxdt\leq\frac{1}{2}\int_{\mathbb{R}^{n}}u_{0}^{2}dx\end{eqnarray}
for every given $T>0$ and $Q_{T}=\mathbb{R}^{n}\times(0,T)$. Let
\begin{eqnarray} G&=&v-u_{\eta}^{m},\nonumber\\
\psi&=&\int_{T}^{t}Gd\tau~~~~~~~~~~~~~~~~~~~~~0<t<T,\nonumber
\end{eqnarray}
where $u_{\eta}$ are the solutions of (1.9). Let $\{
\zeta_{k}\}_{k>1}$ be a smooth cutoff sequence with the following
properties: $ \zeta_{k}(x)\in{}C_{0}^{\infty}(\mathbb{R}^{n})$ and
 \begin{eqnarray}
~~~~~~~~~~~~~~~~~~~~\zeta_{k}(x)=\left\{
\begin{array}{ll}
  1 ~~~~~~~~~~~~~~~~~~~~~~~ ~~~~~~~~~~~~~~~|x|\leq{}k,\\
0< \zeta_{k}(x)<1~~~~~~~~~~~~~~ ~~k<|x|<2k,\\
 0 ~~~~~~~~~~~~~~~~~~~~~ ~~~~~~~~~~~~~~~~~|x|\geq{}2k.
 \end{array}
\right.\nonumber
\end{eqnarray}
    Clearly,  there is a positive constant $ \gamma$ such that
   \begin{eqnarray}   | \nabla{} \zeta_{k}|\leq\frac{ \gamma}{k  }  ~~~~~~\mbox{and}~~~~~~|\Delta{} \zeta_{k}|\leq\frac{ \gamma}{k^{2}}.\end{eqnarray}
Recalling $(v-u_{\eta})_{t}=\Delta{}G$ in $Q_{T}$, we multiply the
equation by $\psi\zeta_{k}$ and integrate by parts
 in $Q_{T}$, we obtain
\begin{eqnarray}
\int_{Q_{T}}(
\zeta_{k}\nabla{}G\cdot\nabla{}\psi{}+\psi\nabla{}G\cdot\nabla{}
\zeta_{k})dxdt= \int_{Q_{T}}(v-u_{\eta})G \zeta_{k}dxdt.
 \end{eqnarray}
Differentiating (3.10) with respect to $T$, we get
\begin{eqnarray}
\int_{\mathbb{R}^{n}}(v-u_{\eta})G
\zeta_{k}dx&=&-\int_{Q_{T}}\left(
\zeta_{k}|\nabla{}G|^{2}+G\nabla{}G\cdot\nabla{}
\zeta_{k}\right)dxdt\nonumber\\
&\leq&-\frac{1}{2}\int_{Q_{T}}\nabla{}G^{2}\cdot\nabla{}
\zeta_{k}dxdt.
 \end{eqnarray}
 Letting $\eta\longrightarrow0$ in (3.11), we get
\begin{eqnarray}
 \int_{\mathbb{R}^{n}}(v-u)(v-u^{m})\zeta_{k}dx
  &\leq&\frac{
\varphi}{k}
\end{eqnarray}
 for some positive $\varphi=\varphi(T)$ thanks to (3.7), (3.8) and (3.9). On the other hand,
\begin{eqnarray}
 \int_{\mathbb{R}^{n}}(v-u)G \zeta_{k}dx
= \int_{\mathbb{R}^{n}}  (v-u )^{2} \zeta_{k}dx+
 \int_{\mathbb{R}^{n}}  (v-u
 )(u-u^{m}) \zeta_{k}dxdt.\nonumber
\end{eqnarray}
Now we   conclude that
\begin{eqnarray}
 \int_{\mathbb{R}^{n}}  (v-u
)^{2} \zeta_{k}dx&\leq&
 \int_{\mathbb{R}^{n}}  |v-u|\cdot
 |u-u^{m}| \zeta_{k}dx
  +\frac{\varphi}{k}\nonumber\\
  &\leq&\frac{1}{2} \int_{\mathbb{R}^{n}} \left[ (v-u
)^{2} + (u-u^{m})^{2}  \right]\zeta_{k}
 dx+\frac{\varphi}{k}.\nonumber
\end{eqnarray}
This implies
\begin{eqnarray}
 \int_{\mathbb{R}^{n}}  (v-u
)^{2} \zeta_{k}dx
  &\leq& \int_{\mathbb{R}^{n}}   (u-u^{m})^{2}  \zeta_{k}
 dx+\frac{2\varphi}{k}\nonumber\\
 &\leq&(m-1)\int_{\mathbb{R}^{n}}  \xi^{m-1} |u-u^{m}|\zeta_{k}
 dx+\frac{2\varphi}{k}\nonumber\\
 &\leq&(m-1)M^{m-1}\int_{\mathbb{R}^{n}}    |u-u^{m}|\zeta_{k}
 dx+\frac{2\varphi}{k}.\nonumber
\end{eqnarray}
Recalling the definition of $\zeta_{k}$, we see $
\int_{\mathbb{R}^{n}}  (v-u )^{2} \zeta_{k}dx\geq
\int_{|x|\leq{}k} (v-u )^{2} dx$. Moreover, $\int_{\mathbb{R}^{n}}
|u-u^{m}|\zeta_{k}dx$ is  bounded thanks to (1.3) and (1.4).
Hence,
  we can get a positive
constant $C_{\ast}=C_{\ast}(T)$ such that
\begin{eqnarray}
 ~~~~~~~~~~~~\int_{|x|\leq{}k} \left[ v(x,t)-u(x,t)
\right]^{2} dx
 &\leq&C_{\ast}\left[(m-1)+\frac{1}{k}\right]\nonumber
\end{eqnarray}
with respect to $t\in(0,T)$ uniformly.~~~~~~~~~~~~~~~~~~~~~~~~~~~~~~~~~~~~~~~~~~~~~~~~~~~~~~~~~~~~~~~~~~~~~~~~~$\square$\\


\end{document}